\documentclass[10pt,twocolumn,twoside]{IEEEtran}

\usepackage[cmex10]{amsmath}
\usepackage{amssymb}
\usepackage[cmex10]{amsmath}
\usepackage{algorithm}
\usepackage{algpseudocode}
\usepackage{graphicx}
\usepackage{color}
\usepackage{subcaption}

\newcommand{\mR}{\ensuremath{\mathbb{R}}}
\newcommand{\Lzi}[1]{\ensuremath{L_{[1,#1)}}}
\newcommand{\Lij}[2]{\ensuremath{L_{(#1,#2)}}}
\newcommand{\Ljn}[1]{\ensuremath{L_{(#1,n]}}}

\newcommand{\Merr}[1]{\ensuremath{R(#1)}}
\newcommand{\Cerr}[1]{\ensuremath{C(#1)}}

\newcommand{\xv}[1]{\ensuremath{\textbf{x}(#1)}}

\newcommand{\xt}[1]{\ensuremath{\overline{x}_{#1}}}
\newcommand{\xtfixed}{\ensuremath{\overline{x}}}
\newcommand{\xtv}{\ensuremath{\overline{\textbf{x}}}}
\newcommand{\Lff}{\ensuremath{L_{\textit{ff}}}}
\newcommand{\Lfl}{\ensuremath{L_{\textit{fl}}}}
\newcommand{\Llf}{\ensuremath{L_{\textit{lf}}}}
\newcommand{\Lll}{\ensuremath{L_{\textit{ll}}}}

\newcommand{\xf}{\ensuremath{\textbf{x}_f(t)}}
\newcommand{\xfd}{\ensuremath{\dot{\textbf{x}}_f(t)}}
\newcommand{\xl}{\ensuremath{\textbf{x}_l(t)}}
\newcommand{\pos}[2]{\ensuremath{p_{#1#2}}}
\newcommand{\posv}{\ensuremath{\textbf{p}}}

\newcommand{\Mmax}{\ensuremath{R_{\text{max}}}}
\newcommand{\Mopt}{\ensuremath{R^*}}
\newcommand{\Ne}[1]{\ensuremath{N(#1)}}
\newcommand{\expec}[1]{\textbf{E}\left[#1\right]}
\newcommand{\tr}[1]{\textbf{tr}\left(#1\right)}

\newcommand{\eqdef}{\ensuremath{\triangleq}}
\newtheorem{theorem}{Theorem}

\newtheorem{problem}{Problem Statement}
\newtheorem{proposition}{Proposition}
\newcommand{\Deg}[1]{\ensuremath{\Delta_{#1}}}

\newcommand{\Gb}{\ensuremath{\overline{G}}}
\newcommand{\Eb}{\ensuremath{\overline{E}}}
\newcommand{\Vb}{\ensuremath{\overline{V}}}
\newcommand{\Wb}{\ensuremath{\overline{W}}}

\newcommand{\Gbb}{\ensuremath{\overline{G}_i}}
\newcommand{\Ebb}{\ensuremath{\overline{E}_i}}
\newcommand{\Vbb}{\ensuremath{\overline{V}_i}}
\newcommand{\Wbb}{\ensuremath{\overline{W}_i}}

\newcommand{\Gbc}{\ensuremath{\overline{G}}}
\newcommand{\Ebc}{\ensuremath{\overline{E}}}
\newcommand{\Vbc}{\ensuremath{\overline{V}}}
\newcommand{\Wbc}{\ensuremath{\overline{W}}}

\newcommand{\Gbcr}{\ensuremath{\overline{G}}}
\newcommand{\Ebcr}{\ensuremath{\overline{E}}}
\newcommand{\Vbcr}{\ensuremath{\overline{V}}}

\newcommand{\Gbbcr}{\ensuremath{\overline{G}_i}}

\newcommand{\Mod}[1]{\ (\text{mod}\ #1)}

\begin{document}
%
\title{Optimal $k$-Leader Selection for Coherence and Convergence Rate in One-Dimensional Networks}
%
%
%

\author{Stacy~Patterson,~\IEEEmembership{Member,~IEEE,}
        Neil~McGlohon,
        and~Kirill~Dyagilev
\thanks{S.~Patterson and N.~McGlohon are with the Department of Computer Science,
Rensselaer Polytechnic Institute, Troy, NY 12180. Email: \emph{sep@cs.rpi.edu}, \emph{mcglon@rpi.edu}}
\thanks{K.~Dyagilev is with the Department of Computer Science, Johns Hopkins University, Baltimore, MD 21218. Email: \emph{kirilld@cs.jhu.edu}}}

\maketitle


\begin{abstract}
We study the problem of optimal leader selection in consensus networks under two performance measures (1) formation coherence when subject to additive perturbations, as quantified by the steady-state variance of the deviation from the desired trajectory, and (2) convergence rate to a consensus value. The objective is to identify the set of $k$ leaders that optimizes the chosen performance measure. In both cases,  an optimal leader set can be found by an exhaustive search over all possible leader sets; however, this approach is not scalable to large networks. In recent years, several works have proposed approximation algorithms to the $k$-leader selection problem, yet the question of whether there exists an efficient, non-combinatorial method to identify the optimal leader set remains open. This work takes a first step towards answering this question. We show that, in one-dimensional weighted graphs, namely path graphs and ring graphs, the $k$-leader selection problem can be solved in polynomial time (in both $k$ and the network size $n$). We give an $O(n^3)$ solution for optimal $k$-leader selection in path graphs and an $O(kn^3)$ solution for optimal $k$-leader selection in ring graphs.
\end{abstract}



\section{Introduction}
\IEEEPARstart{W}{e} explore the problem of leader selection in leader-follower consensus systems.  Such systems arise in the context of
vehicle formation control~\cite{RBM05}, distributed clock synchronization~\cite{EKPS04}, and distributed localization in sensor networks~\cite{BH09}, among others.
In these systems, several agents act as \emph{leaders} whose state trajectory serves as the reference for the entire system.  These leaders may be controlled autonomously or by a system owner.
The remaining agents are \emph{followers}.  Each follower updates its state based on relative measurements of the states of its neighbors.
The objective of the leader-follower system is for the entire formation to maintain a desired global state.

We consider leader selection under two different dynamics.  In the first, which we call \emph{noisy formation dynamics},
the follower agents' measurements are corrupted by stochastic noise.  With these additive perturbations, 
the agents cannot maintain the formation exactly, however,
the steady-state variance of the deviation of the agents' states from the desired states is bounded~\cite{BH06}.  This variance is related to the \emph{coherence} of the formation and is quantified by an $H_2$-norm of the leader-follower system~\cite{PB10,B12}.  
It has been shown that the coherence depends  on which agents act as leaders~\cite{BH06,PB10}, and so, by judiciously choosing the leader set, one can minimize the total variance of the formation.
   The \emph{$k$-leader selection problem for coherence} is precisely to select a set of at most $k$ leaders that minimize this total variance.   The second dynamics we consider are \emph{consensus dynamics}.  In this case, the leaders all share the same target value, and the followers' relative measurements are exact.  Under these dynamics, the network topology and the choice of the leader set  determine the convergence rate.
 The \emph{$k$-leader selection problem for fast convergence} is to, given a consensus network, select a set of at most $k$ leaders so that the convergence rate is maximized.  Both leader selection problems can be solved by an exhaustive search over all subsets of agents of size at most $k$, but this approach is not tractable in large networks for anything other than small values of $k$.

The leader selection problem for coherence has received a great deal of attention in recent years.  
In~\cite{CBP14}, the authors show that the total variance of the deviation from the desired trajectory is a super-modular set function~\cite{N78}.
 This super-modularity property implies that one can use a greedy, polynomial-time algorithm to find a leader set for which the coherence is within a provable bound of optimal.
 In~\cite{LFJ14}, the authors give algorithms that yield lower and upper bounds on coherence of the optimal leader set.
 Another recent work~\cite{FL13} has shown a connection between the optimal leader set for coherence and information centrality measures.  They have used this connection to give efficient algorithms for finding
the optimal single leader and optimal pair of leaders in several network topologies.  

With respect to the convergence rate in leader-follower consensus networks, 
recent works have characterized this convergence rate in terms of the spectrum of a 
weighted Laplacian matrix~\cite{PMB08,RJME09,GS13}.  Several works have developed bounds for the convergence rate based on the network topology~\cite{GS13,PS14,CABP14}.  Finally, the leader selection problem for fast convergence has been studied in~\cite{CABP12,CABP14}.  
In these works, the authors propose relaxations of the leader selection problem that can be solved efficiently.  While these relaxed formulations perform well in evaluations, there are no guarantees on the optimality of their solutions.

Despite this recent interest in leader selection problems, the question of whether it is possible to find an \emph{optimal} leader set for coherence or convergence rate using a non-combinatorial approach remains open. 
In this work, we take a first step towards answering this question.  
We show that, in one-dimensional weighted graphs, specifically, path graphs and ring graphs, the $k$-leader selection problems for coherence and fast convergence can be solved in time that is polynomial in both $k$ and the network size $n$.  For the coherence problem, we transform the leader selection problem into a problem of finding a minimum-weight path
in a weighted directed graph, a problem that can be solved in polynomial time.
We then use a slightly modified version of the well-known Bellman-Ford algorithm~\cite{CLRS10}
to find this minimum-weight path and thus the optimal leader set.
For the fast convergence problem, we transform the problem into a task of finding the widest path in a weighted directed graph. 
A version of the Bellman-Ford algorithm can then be used to find the optimal leader set for fast convergence in polynomial time.
 For path graphs, our algorithms find the optimal leader set of size at most $k$ in
 $O(n^3)$ time.  In ring graphs, our algorithms find the optimal leader set of size at most $k$ in $O(kn^3)$ time.
Our algorithm for optimal $k$-leader selection for coherence first appeared in~\cite{PMD14}.  
This paper extends our previous work to the leader selection problem for fast convergence.

The proposed algorithms have applications in one dimensional leader-follower networks such as perimeter patrolling, where a ring of autonomous agents 
encircles a hazard such as a fire or an oil spill to monitor its location and conditions~\cite{susca2014synchronization,BaseggioCMPS10}.
In these applications, a system operator may issue commands to leaders to change the velocity of the patrol or the distance from the hazard, for example.
Other relevant applications include   on-road and underwater one-dimensional vehicular platoons~\cite{wang2012coordinated,1431191}, where leader vehicles
may dictate the platoon trajectory.

Our work was inspired by recent work on network facility location~\cite{HD02}.
In network facility location, the objective is to identify a subset of nodes to
serve as facilities that minimizes a function of the network distances between the remaining nodes and their closest facility nodes.
In leader selection, the performance of the network for a given leader set
also depends on the locations of the follower nodes with respect to the leaders, though in a different way than in classical facility location.
Our  algorithms use a similar approach to the facility location algorithm for a real line that was presented in~\cite{CW14},
where the authors solve the facility location problem by reducing it to a minimum-weight path problem over a directed acyclic graph.
We note that our graph construction and our path-finding algorithm are both significantly different from those in~\cite{CW14}.  Nevertheless, this shows that
it is possible to exploit connections between facility location and leader selection. We anticipate that
this connection can be used to develop efficient leader selection algorithms for other network topologies that have efficient facility location algorithms,
for example, tree graphs~\cite{KH79}.

The remainder of this paper is organized as follows.  In Section~\ref{problem.sec}, we present the system model and $k$-leader selection problems.
In Section~\ref{algnoise.sec}, we present our  algorithms for finding the optimal leader set for coherence in path and ring graphs.
In Section~\ref{algconsensus.sec}, we present our algorithms for finding the optimal leader set for fast convergence in path and ring graphs.
Section~\ref{examples.sec} gives numerical examples comparing the optimal leader set to the leader set selected by a greedy algorithm.
Finally, we conclude in Section~\ref{conclusion.sec}.

\section{Problem Formulation} \label{problem.sec}

In this section, we describe the dynamics of the leader-follower systems and formally define the $k$-leader selection problems.

\subsection{System Model}
We consider a network of agents, modeled by a directed, connected  graph $G=(V,E)$, where the node set $V = \{1, 2, \ldots, n\}$ represents the agents.
We use the terms node and agent interchangeably.  The edge set $E$ represents the communication structure.
We assume that if $(i,j) \in E$, then $(j,i) \in E$; in this case, both $i$ and $j$ can exchange information.
We denote the neighbor set of a node $i$ by $\Ne{i}$, i.e., $\Ne{i} \eqdef \{ j \in V~|~(i,j) \in E\}$.
In this work, we restrict our study to one-dimensional graphs,
namely path graphs and ring graphs.
Without loss of generality, for path graphs, we assume the node IDs are assigned in order along the path, and 
for ring graphs, we assume that the node IDs are assigned in ascending order around the ring in a clockwise fashion.

Each agent $i$ has a state $x_i(t) \in \mR$.  The objective is for each pair of neighboring agents $i$ and $j$ to maintain a pre-specified difference $\pos{i}{j}$ between their states,
\begin{equation} \label{diff.eq}
x_i(t) - x_j(t) = \pos{i}{j}~~~~\text{for each}~(i,j) \in E.
\end{equation}
If $x_i(t)$ represents an agent's position, for example, $\pos{i}{j}$ is the desired distance between agents $i$ and $j$.
We assume that each agent $i$ knows $\pos{i}{j}$ for all $j \in \Ne{i}$.

A subset $S \subseteq V$ of agents act as leaders.  The  states of these agents
serve as reference states for the network.   The state of each leader $s \in S$ remains fixed at its reference value
$\xt{s}$.
The remaining agents $v \in V \setminus S$ are followers. A follower  updates its state based on measurements of the differences between its state and the states
of its neighbors.  Let $\xtv$ denote the vector of states that satisfy  (\ref{diff.eq}) when the leader states are fixed at their reference values.
We assume that at least one such $\xtv$ exists.

\subsubsection{Noisy formation dynamics}
We first consider the formation dynamics presented in~\cite{CBP14}, which we call \emph{noisy formation dynamics}.
In these dynamics, the followers' measurements are corrupted by stochastic noise.
The dynamics of each follower agent is given by,
\begin{equation}\label{noise.eq}
\dot{x}_i(t) = -\sum_{j \in \Ne{i}} W_{ij}  \left(x_i(t) - x_j(t) - \pos{i}{j}+ \epsilon_{ij}(t)\right),
\end{equation}
where $W_{ij}$ is the weight for link $(i,j)$, and $\epsilon_{ij}(t)$ are  zero-mean white noise processes with autocorrelation functions
$\expec{\epsilon_{ij}(t)\epsilon_{ij}(t + \tau)} = \nu_{ij} \delta(\tau)$.  Here $\delta(\cdot)$ denotes the unit impulse function.
 Each $\epsilon_{ij}$ is independent, and
$\epsilon_{ij}$  and $\epsilon_{ji}$ are identically distributed.
As specified in~\cite{CBP14}, the link weights are given by $W_{ij} = \frac{\nu_{ij}}{\Deg{i}}$, where $\Deg{i} = \sum_{j \in \Ne{i}}(1/\nu_{ij})$.
This policy ensures that the steady-state variance of the deviation of the follower states from a desired $\xtv$ is bounded~\cite{CBP14}.
We discuss this in more detail in Section~\ref{formationdyn.sec}.

Let $L$ be a weighted Laplacian matrix of the graph $G$, with each component defined as,
\begin{equation} \label{formationL.eq}
L_{ij} = \left\{ \begin{array}{ll}
- \frac{1}{\nu_{ij}} & ~\text{if}~(i,j) \in E \\
\Deg{i} &~\text{if}~i=j \\
0 &~\text{otherwise.}
\end{array} \right.
\end{equation}
Arranging the node states $\xv{t}$ so that $\xv{t} =  [ \xf^T~\xl^T]^T$,
where $\xf$ contains the follower nodes' states and
$\xl$ contains the leader nodes' states,
$L$ can be decomposed as,
\begin{equation} \label{Lff.eq}
L = \left[ \begin{array}{cc}
\Lff & \Lfl \\
\Llf & \Lll
\end{array} \right].
\end{equation}
The sub-matrix $\Lff$ defines the interactions between followers, and the sub-matrix $\Lfl$ defines the impact of the leaders on the followers.
Since $G$ is connected, $\Lff$ is positive definite~\cite{mesbahi2010graph}.

Let $B$ be the $|V - S| \times |E|$  matrix where  each column $k$ corresponds to an
edge $(i,j) \in E$.  $B_{ik} = L_{ij}$ if edge $k$ is the edge $(i,j)$, and $B_{ik} = 0$ otherwise.
With these definitions, the dynamics of the follower agents can be written as,
\[
\xfd = -D_f^{-1}(\Lff \xf + \Lfl \xl + B\posv ) + \epsilon(t),
\]
where $D_f$ is a $|V - S| \times |V - S|$ diagonal matrix with diagonal entries $\Deg{i}$, $\posv$ is the vector of desired differences $\pos{i}{j}$,
and $\epsilon(t)$ is the $|V-S|$-vector of noise processes.
Let $\xtv_l$ be the vector of reference states of the leader nodes.
The desired state of the follower nodes is given by $\xtv_f = - \Lff^{-1}(\Lfl \xtv_l + B\posv)$.
We refer the reader to~\cite{CBP14} for the details of this derivation.

\subsubsection{Consensus dynamics}
We also consider a simplified leader-follower system where $\pos{i}{j} = 0$ for all $(i,j) \in E$ and 
where $\xtv_{l} = \xtfixed \textbf{1}$ with $\xtfixed \in \mR$, i.e., all leaders share the same reference
 value $\xtfixed$.  
Note that $\xtfixed$ may be non-zero.
 Each follower updates its state based on exact relative measurements of the states of its neighbors, 
 \begin{equation} \label{noiseless.eq}
 \dot{x}_i(t) = -\sum_{j \in \Ne{i}} W_{ij}  \left(x_i(t) - x_j(t) \right).
 \end{equation}
We assume that for all follower nodes $i$, $W_{ij} > 0$, and that for all follower nodes  $i$ and $j$, $W_{ij} = W_{ji}$.
In the absence of leader nodes, these dynamics correspond to the standard consensus dynamics over an undirected graph.
We  thus call the dynamics \emph{consensus dynamics}. 
The objective in this setting is for all follower agents to converge to the leaders' value $\xtfixed$.

We define the weighted Laplacian matrix $L$ as
\begin{equation} \label{consensusL.eq}
L_{ij} = \left\{ \begin{array}{ll}
- W_{ij} & ~\text{if}~(i,j) \in E \\
\sum_{j \in N(i)} W_{ij} &~\text{if}~i=j \\
0 &~\text{otherwise.}
\end{array} \right.
\end{equation}
Let $\Lff$ be the principle sub-matrix of $L$ capturing the follower-follower interactions, as in (\ref{Lff.eq}).
The followers' states evolve as,
\begin{equation} \label{cnsnfollow.eq}
\dot{\textbf{x}}_f(t) = -\Lff \xf - \Lfl \xtv_{l},
\end{equation}
Defining $\tilde{\textbf{x}}_f(t)  = \xf - \xtv_{f}$,  where $\xtv_{f} = \xtfixed \textbf{1}$, we rewrite (\ref{cnsnfollow.eq}) as,
\begin{align*}
\textstyle \frac{d}{dt} \left( \tilde{\textbf{x}}_f(t)  + \xtv_f \right) &= -\Lff  \left( \tilde{\textbf{x}}_f(t)  + \xtv_f \right) - \Lfl \xtv_{l} \\
\textstyle \frac{d}{dt} \tilde{\textbf{x}}_f(t) &= -\Lff \tilde{\textbf{x}}_f(t).
\end{align*}
where the second equality follows from the fact that $\Lff \textbf{1} + \Lfl \textbf{1} = 0$. The rate of convergence to consensus
is thus equivalent to the rate at which  $\tilde{\textbf{x}}_f(t)$ converges to 0~\cite{PB10}.

\subsection{Leader Selection Problems}

\subsubsection{Leader Selection with Noisy Formation Dynamics} \label{formationdyn.sec}

Under the  dynamics (\ref{noise.eq}), the followers' states deviate from $\xtv$, however, the variances of these deviations are bounded in the mean-square sense.
For a follower agent $i$, let $r_i$ be the steady-state variance of the deviation from $\xt{i}$,
\[
r_i \eqdef \lim_{t \rightarrow \infty} \expec{(x_i(t) - \xt{i})^2}.
\]
It has been shown that,
\begin{equation} \label{ri.eq}
r_i = \textstyle \frac{1}{2} (\Lff^{-1})_{ii},
\end{equation}
where $\Lff$ is as defined in (\ref{Lff.eq}), i.e., $\Lff$ is the sub-matrix of the Laplacian $L$ defined in~(\ref{formationL.eq}) with the rows and columns corresponding
to nodes in $S$ removed.  A derivation of (\ref{ri.eq}) can be found in~\cite{BH06}.

We measure the performance of the formation for a leader set $S$ by the total steady-state variance of the deviation from $\xtv$,
\begin{equation} \label{Merr.eq}
\Merr{S} \eqdef \textstyle \sum_{i \in V \setminus S} r_i = \textstyle \frac{1}{2} \tr{\Lff^{-1}}.
\end{equation}
A formation with a small $\Merr{S}$ exhibits good coherence, meaning the formation closely resembles a rigid formation.

A natural question that arises is how to identify a leader set $S$ of a certain size that minimizes the total steady-state variance  (\ref{Merr.eq}).
This question is formalized as the \emph{$k$-leader selection problem for coherence}.
\begin{problem} \label{coherence.prob}
The \emph{$k$-leader selection problem for coherence} is
\begin{equation} \label{problem.eq}
\begin{array}{ll}
\text{minimize}~& \Merr{S} \\
\text{subject to}~& |S|  \leq k.
\end{array}
\end{equation}
\end{problem}

A na\"{\i}ve solution to this problem is to construct all subsets of $V$ of size at most $k$, evaluate $R(S)$ for each leader set,
and choose the set for which $R(S)$ is minimized.  The computational complexity of this solution is combinatorial since the number of leader sets that would need to be evaluated is
$\sum_{i=1}^k {n \choose i}$.  In Section~\ref{algnoise.sec}, we show that for one-dimensional graphs, the optimal leader set can be found with time complexity that is polynomial in both $n$ and $k$.

\subsubsection{Leader Selection with Consensus Dynamics}

Under the leader-follower consensus dynamics (\ref{noiseless.eq}), when all leader states are equal, the follower states converge asymptotically to $\xtfixed$~\cite{RJME09}.
We quantify the performance of this system in terms of the convergence rate.
The convergence rate is determined by the smallest eigenvalue of $\Lff$~\cite{horn2012matrix}.
This eigenvalue, in turn, depends on which agents are selected as leaders.  We let $\Cerr{S}$ denote the smallest eigenvalue of $\Lff$ for the
leader set $S$.  The \emph{$k$-leader selection problem for fast convergence} is to find the set $S$, of size at most $k$, that maximizes this eigenvalue.
\begin{problem}
The \emph{$k$-leader selection problem for fast convergence} is
\begin{equation} \label{convergence.prob}
\begin{array}{ll}
\text{maximize}~& \Cerr{S} \\
\text{subject to}~& |S|  \leq k.
\end{array}
\end{equation}
\end{problem}

As with Problem~\ref{coherence.prob}, this problem can be solved by a search over all possible leader sets of size $k$ or less, an approach which has combinatorial complexity.  In Section~\ref{algconsensus.sec}, we give polynomial-time algorithms that find the optimal leader set in one-dimensional graphs.

\begin{figure}
\begin{center}
\includegraphics[scale=.35]{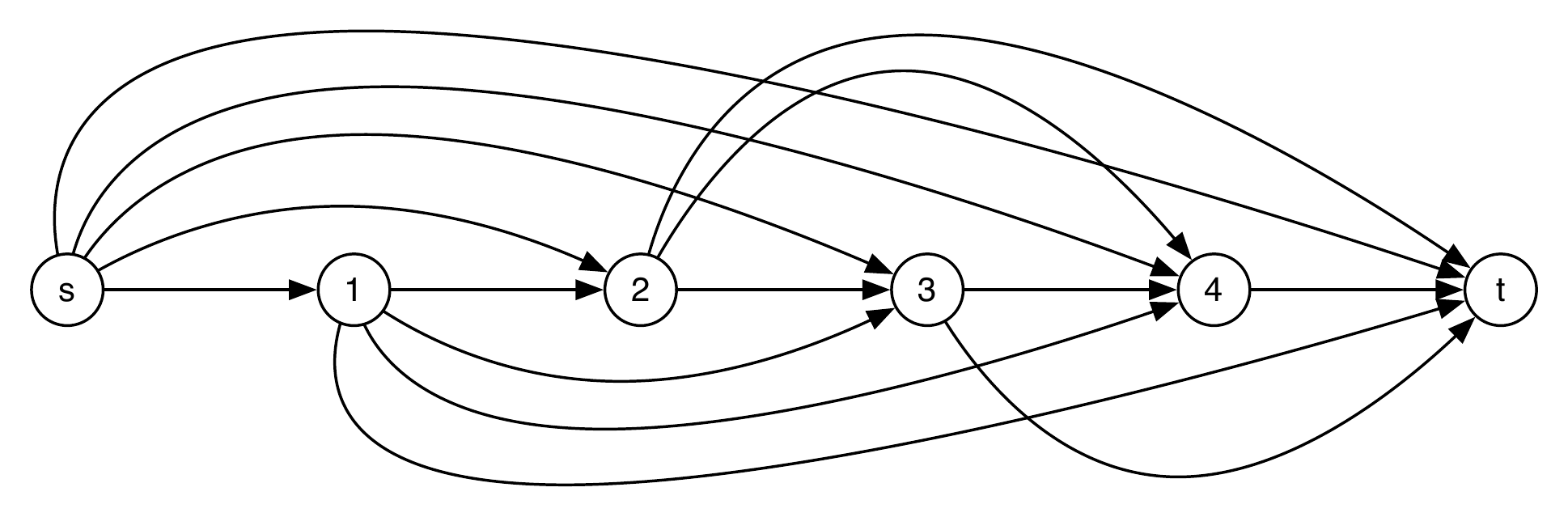}
\end{center}
\caption{Digraph generated from four node path graph.}  \label{dag.fig}
\end{figure}

\section{Optimal Leader Selection for Coherence} \label{algnoise.sec}

In this section, we present our leader selection algorithms for Problem~\ref{coherence.prob} for path and ring graphs.
Our approach for both graph types is to first reduce the leader selection problem to the task of finding a minimum-weight path in a
digraph. 
We then solve this minimum-weight path problem using a modified version of the
Bellman-Ford algorithm~\cite{CLRS10}.

\subsection{Optimal $k$-Leader Selection for a Path Graph} \label{algnoisepath.sec}

In a path graph with $k$ leader nodes, the matrix $\Lff$ is block diagonal with at most $k+1$ blocks.  Further,
each block is tridiagonal.
For example, consider a path graph with leader set $S = \{ \ell_1, \ell_2, \ldots, \ell_k\}$, where $\ell_i < \ell_{i+1}$ for $i=1, \ldots, (k-1)$.
The matrix $\Lff$ is formed from  the weighted Laplacian $L$ by removing the rows and columns corresponding to the nodes in $S$.
$\Lff$ can be written as,
\begin{equation} \label{blockL.eq}
\Lff = \left[  \begin{array}{ccccc}
\Lzi{\ell_1} & \textbf{0} & \cdots & \cdots & \textbf{0} \\
\textbf{0} & \Lij{\ell_1}{\ell_2} & \textbf{0} & \cdots & \vdots \\
\vdots &  \ddots & \ddots & \ddots &   \vdots \\
\textbf{0} & \ddots & \ddots & \Lij{\ell_{k-1}}{\ell_k} & \textbf{0} \\
\textbf{0} & \cdots & \cdots & \textbf{0} & \Ljn{\ell_k}
\end{array} \right],
\end{equation}
where the matrices $\Lzi{\ell_1}$, $\Lij{\ell_i}{\ell_{i+1}}$, $i=1, \ldots ,(k-1)$, and $\Ljn{\ell_k}$ are defined as follows.
$\Lzi{\ell_1}$ is the sub-matrix of $L$ consisting of the rows and columns $1$ through $\ell_1-1$ of $L$.
This matrix models the interactions of the follower nodes $i =1, \ldots, (\ell_1-1)$.
Note that node $\ell_1$ is the only leader that affects the states of these nodes.  
$\Lij{\ell_i}{\ell_{i+1}}$  is the sub-matrix of $L$ consisting of the rows and columns indexed from $\ell_i+1$ to $\ell_{i+1}-1$, inclusive.
This matrix models the interactions of the nodes between leader node $\ell_i$ and leader node $\ell_{i+1}$;
these follower nodes are influenced by both leaders.
Finally, $\Ljn{\ell_k}$ consists of the rows and columns $i=(\ell_k+1),\ldots, n$ of $L$, inclusive, and models the interactions of the follower nodes $i=(\ell_k+1), \ldots, n$.
Leader node $\ell_k$ is the only leader that influences the states of these nodes.
If there are leaders $u, v \in S$ such that $v = u+1$, then the corresponding sub-matrix of $L$ will be of size 0.

We  construct a weighted digraph $\Gb = (\Vb, \Eb, \Wb)$ based on $\Lff$ as follows.
The set of nodes consists of $V$ and an additional source node $s$ and target node $t$, i.e.,
$\Vb = V \cup \{s,t\}$, as shown in Figure~\ref{dag.fig}.
The edge set $\Eb$  contains edges from $s$ to every node $v \in V$.  The weight of edge $(s,v)$ is
$\frac{1}{2}\tr{{\Lzi{v}}^{-1}}$.  This edge weight is the total steady-state variance  for nodes $1, \ldots, (v-1)$ when
node $v$ is a leader node and there are no other leader nodes $u$ with $u < v$.  
$\Eb$ also contains edges from each node $u \in V$ to each node $v \in V$ with $u < v$.  
The weight of edge $(u,v)$ is $\frac{1}{2}\tr{{\Lij{u}{v}}^{-1}}$.  This weight is the total variance of the nodes
$i=(u+1), \ldots,(v-1)$ when nodes $u$ and $v$ are leaders and there are no other leader nodes $w$ with $u < w < v$.
Finally,  $\Eb$ also contains edges from every node $v \in V$ to node $t$.  The weight of edge $(v,t)$ is $\frac{1}{2}\tr{{\Ljn{v}}^{-1}}$.
This weight is the total variance of the nodes $i={v+1}, \ldots, n$ when $v$ is a leader and there are no other leader nodes $u$ with
$u > v$.  The weights of edges $(s,1)$, $(n,t)$, and $(u,u+1),~u=1 \ldots (n-1)$ are 0.

\begin{proposition} \label{path.prop}
Let $G$ be a path graph, and let $L$ be the respective weighted Laplacian defined in (\ref{formationL.eq}).
Let $\Gb$ be the digraph generated from $G$ and $L$.
Further, let $P=\{s, u_1, u_2, \ldots, u_k, t\}$ be the set of nodes on a path from $s$ to $t$ in $\Gb$ that contains $k+1$ edges, and let $w$ be the corresponding path weight.
Then, for $S  = P \setminus \{s,t\}$, the coherence for the leader set $S$ in the graph $G$ is $R(S) = w$.
\end{proposition}

\begin{IEEEproof}
For $S = P \setminus \{s,t\}$, $R(S) = \frac{1}{2} \tr{\Lff^{-1}}$, where $\Lff$ is the submatrix of the Laplacian where the rows and columns of nodes in $S$ have been removed.
Using the definition of $\Lff$ in (\ref{blockL.eq}), we can write $R(S)$ as
\begin{align} \label{R1d.eq}
R(S) &= \textstyle \frac{1}{2} \tr{{\Lzi{u_1}}^{-1}} + \frac{1}{2} \sum_{i=1}^{k-1} \tr{ {\Lij{u_i}{u_{i+1}}}^{-1}}  \nonumber  \\
&~~~~~+ \textstyle \frac{1}{2} \tr{{\Ljn{u_k}}^{-1}}.
\end{align}
Note that if any sub-matrix of $L$ is of size 0, then the trace of the inverse of this sub-matrix is 0.  

In the digraph $\Gb$, the weight of the path formed by $P$ is 
\[
w = w_{s,u_1} + \sum_{i=1}^{k-1} w_{u_i,u_{i+1}} + w_{u_k,t},
\]  
where $w_{i,j}$ is the weight of the edge from $i$ to $j$ in $\Gb$.
The edge weights in the digraph are defined so that $w_{s,u_1} = \frac{1}{2} \tr{{\Lzi{u_1}}^{-1}}$,
$w_{u_i, u_{i+1}} = \frac{1}{2} \tr{ {\Lij{u_i}{u_{i+1}}}^{-1}} $ for $i=1,\ldots,k+1$, and $w_{u_k,t} = \frac{1}{2} \tr{{\Ljn{u_k}}^{-1}}$.
Therefore, $R(S) = w$.
\end{IEEEproof}

It follows from Proposition~\ref{path.prop} that, to find a leader set $S$,  with $|S| \leq k$, that minimizes (\ref{R1d.eq}),  one seeks a minimum-weight path from $s$ to $t$ in
$\Gb$ that has at most $k+1$ edges.  The optimal leaders  are the nodes along this path between $s$ and $t$.  To find this minimum-weight path,
we use a slightly modified implementation of the Bellman-Ford algorithm~\cite{CLRS10}.
The Bellman-Ford algorithm is an iterative algorithm that finds the minimum-weight paths (of any length) from a source node to every other node in the graph.
While there are more efficient algorithms that solve this same problem, Bellman-Ford offers the benefit that, in each iteration $m$, the algorithm finds the minimum-weight paths of $m$ edges.
Therefore, we can execute the Bellman-Ford algorithm for $k+1$ iterations to find the minimum-weight path of at most $k+1$ edges.
We have made slight modifications to this algorithm to make it possible to retrieve not only the weight of the path but the list of nodes traversed in this path.
Our modified version of the Bellman-Ford algorithm is detailed in the appendix.

\begin{algorithm}[t]
\small
\caption{Algorithm for finding the optimal solution to the $k$-leader selection problem for coherence in a path graph.} \label{optpath.alg}
\begin{algorithmic}[1] \footnotesize
\State \textbf{Input:} $G = (V,E)$, edge weights $W_{ij}$, $k$
\State \textbf{Output:} Set of leader nodes $S$, error $R(S)$
\Statex ~
\State $L \gets$ Laplacian of $G$
\State $\Gb \gets $ digraph constructed from $G$ and $L$
\State {\emph{/*~Get node set and weight of min. weight path*/}}
\State $(P, {weight}) \gets \textsc{ModifiedBellmanFord}(\Gb, s, t, k+1)$ \label{bf.line}
\Statex ~
\State \emph{/*~Construct leader set from min. weight path.~*/}
\State $S \gets P \setminus \{s,t\}$
\State $R(S) \gets {weight}$
\State \textbf{return} $(S, R(S))$
\end{algorithmic}
\end{algorithm}

The pseudocode for our $k$-leader selection algorithm for coherence is given in Algorithm~\ref{optpath.alg}.  The algorithm returns the optimal set of leaders of size at most $k$.
A leader set may have cardinality $h < k$ if the inclusion of more than $h$ leaders does not decrease $R(S)$.

\begin{theorem}
For a path graph $G$ with $n$ nodes, the $k$-leader selection algorithm identifies the leader set $S$, with $|S| \leq k$, that minimizes $\Merr{S}$
in $O(n^3)$ operations.
\end{theorem}
\begin{IEEEproof}
From Proposition~\ref{path.prop} and the correctness of the Bellman-Ford algorithm, Algorithm~\ref{optpath.alg} finds the optimal leader set.
We next examine the complexity of the proposed algorithm.

The algorithm consists of two phases.  The first phase is the construction of the digraph $\Gb = (\Vb, \Eb)$.
The edge set $\Eb$ consists of $n$ edges with $s$ as their source (one to each $v \in V$), $n$ edges with $t$ as their sink (one from each $v \in V$),
and one edge from each $u \in V$ to each $v \in V$ with $u < v$.  Thus $| \Eb | \in O(n^2)$.
To find each edge weight, we must find the diagonal entries of the inverse of a tridiagonal matrix of size at most $(n-1) \times (n-1)$.
These diagonal entries can be found in $O(n)$ operations~\cite{RH91}.  Therefore the digraph $\Gb$ can be constructed in $O(n^3)$ operations.

The second phase of the algorithm is to find the shortest path of at most $k+1$ edges from $s$ to $t$ in $\Gb$.  For a graph with $m$ edges,
the Bellman-Ford algorithm finds the minimum-weight-path of length at most $h$  edges in $O(hm)$ operations~\cite{CLRS10}.
Therefore, the second phase of our algorithm has time complexity $O(kn^2)$.

Combining the two phases of the algorithm we arrive a total time complexity of $O(n^3)$.
\end{IEEEproof}

\begin{algorithm}[t]
\small
\caption{Algorithm for finding the optimal solution to the $k$-leader selection problem for coherence in a ring graph.}  \label{optring.alg}
\begin{algorithmic}[1] \footnotesize
\State \textbf{Input:} $G = (V,E)$, edge weights $W_{ij}$, $k$
\State \textbf{Output:} Set of leader nodes $S$, error $R(S)$
\Statex~
\State $L \gets $ weighted Laplacian for $G$

\State $minWeight \gets \infty$
\State $P \gets \bot$
\For{$i=1 \ldots n$}
\State  $\Gbb \gets$ digraph constructed from $G$ and $L$ for candidate leader $i$
\State {\emph{/*~Get node set and weight of min. weight path}}
\State {\emph{~~~~~~~~~that contains node $i$*/}}
\State $(P, weight) \gets \textsc{ModifiedBellmanFord}(\Gbb, s_i, t_i, k-1)$
\If{$weight < minWeight$}
	\State $minWeight \gets weight$
	\State $minP \gets P$
\EndIf
\EndFor
\Statex~
\State \emph{/*~Construct leader set from min. weight path.~*/}
\State $S \gets minP \setminus \{s_i, t_i \}$
\State {$R(S) \gets minWeight$}
\State \textbf{return} $(S, R(S))$
\end{algorithmic}
\end{algorithm}

\subsection{Optimal $k$-Leader Selection for a Ring Graph} \label{algnoisering.sec}
In a ring graph with $k$ leaders, the leader-follower system can be decomposed into $k$ independent subsystems (with some possibly consisting of zero nodes).
Each of these subsystems corresponds to a segment of the graph where two leader nodes $u$ and $v$ form the boundaries of this segment,
 where $v$ follows $u$ in the clockwise direction, and where there are no other leader nodes in this segment.  Nodes $u$ and $v$ are not included in the subsystem.
The coherence of the
subsystem is given by the sub-matrix of the Laplacian consisting of the rows and columns corresponding to the nodes between $u$ and $v$ in the ring.
We denote this submatrix by $L_{(u \to v)}$.
We further let $R_{(u\to v)}$ denote the the total steady-state variance for this subsystem, 
$R_{(u \to v)} = \frac{1}{2} \tr{L_{(u \to v)}^{-1}}$. For example, in a ring graph with $n>6$ nodes, the matrix $L_{(5 \to 2)\}}$ contains rows and columns that correspond to nodes $6, 7, \ldots, n, 1$. 
The value $R_{(5\to 2)}$ is equal to the total steady-state variance of the nodes $6,7,...,n,1$ when nodes $5$ and $2$ are leaders and nodes $\{6,7,\ldots,n,1\}$ are followers.

To find the optimal leader set of size at most $k$, we first select one node $i$ as a candidate leader.
We then translate the problem finding the remaining $k-1$ leaders into a problem of finding a minimum weight path of at most $k$ edges
over a weighted digraph.
The digraph is described below.
To ensure that our algorithm finds the optimal leader set, the algorithm performs this translation and path-finding for each possible initial leader $i= 1, \ldots, n$.
The optimal leader set is the set with the minimum weight path among these $n$ minimum weight paths (one for each initial leader selection).
Pseudocode for the algorithm is given in Algorithm~\ref{optring.alg}.

For a given initial candidate leader $i$, its weighted digraph
$\Gbb = (\Vbb, \Ebb, \Wbb)$ is defined as follows.  The vertex set of $\Vbb$ contains a source node $s_i$, a target node $t_i$,
and the vertices in $V$ excepting $i$, i.e.,
$\Vbb = \{s_i, t_i\} \cup ( V \setminus \{i\})$.
The edge set $\Ebb$  contains directed edges from $s_i$ to every node $v \in ( V \setminus \{i\})$.
The weight of an edge $(s_i,v)$ is $w_{s_i,v} =  R_{(i \to v)}$.
$\Ebb$ also contains edges from from each node $u,v \in ( V \setminus \{i\})$.
The weight of an edge $(u,v)$ is $w_{u,v} = R_{(u\to v)}$.
Finally, $\Ebb$  contains directed edges from every node $v \in (V \setminus \{i\})$ to $t_i$
with weights $w_{v,t_i} = R_{(v \to i)}$.

\begin{proposition} \label{ring.prop}
Let $G$ be a ring graph, and let $L$ be the respective Laplacian as defined in (\ref{formationL.eq}).
Let $\Gbb$ be the weighted digraph generated from $G$ and $L$ for a given node $i$.
Further, let $P=\{s_i, u_1, u_2, \ldots, u_{k-1}, t_i\}$ be the set of nodes on a path from $s_i$ to $t_i$ in $\Gbb$ that contains $k$ edges, and let $w$ be the corresponding path weight.
Then, for $S  = P \setminus \{s_i,t_i\} \cup \{i\}$, $\Merr{S} = w$.
\end{proposition}

\begin{IEEEproof}
For a leader set $S$ that contains $i$,
\[
R(S) = R_{(i \to u_1)} + \sum_{j=1}^{k-2} R_{(u_j \to u_{j+1})}  + R_{(u_{k-1} \to i)}.
\]
This is precisely $w$, the weight of path $P$ in $\Gbb$.
\end{IEEEproof}

It follows from Proposition~\ref{ring.prop} that, to find the optimal leader set that contains node $i$, one must find the minimum weight path from
$s_i$ to $t_i$ in $\Gbb$ of at most $k$ edges.  The weight of this path is the minimal $\Merr{S_i}$ when $i \in S_i$ and $|S_i| \leq k$.
Our algorithm uses the modified Bellman-Ford algorithm to find this path for each $i$.
Let $S_i$ be the set of nodes along this path.
The optimal leader set for the graph is then
\[
S^* = \arg \min_{S_i, i \in V} \Merr{S_i}.
\]

\begin{theorem} \label{coherence.thm}
For a ring graph $G$ with $n$ nodes, the  $k$-leader selection algorithm identifies the leader set $S$, with $|S| \leq k$, that minimizes $\Merr{S}$
in $O(kn^3)$ operations.
\end{theorem}
\begin{IEEEproof}
The fact that the leader set $S$ is optimal follows from Proposition~\ref{ring.prop} and the correctness of the Bellman-Ford algorithm.

With respect to the computational complexity, the algorithm constructs $n$ weighted digraphs $\Gbb$, $i=1,\ldots, n$, and a shortest-path algorithm is executed on each digraph.

To construct these digraphs, first, the weight of each edge $(u,v)$, $u,v \in \Vbb, u \neq v$ is computed.
To compute the weight for edge $(u,v)$, where $v$ follows $u$ on the ring in the clockwise direction,
 we first construct the matrix $M$ from $L$ by shifting the rows and columns of $L$ so that node
$u$ corresponds to the first row and column of $M$.
The index (row and column of $M$) corresponding to a node $u$ after this shift is $1$,
and the index corresponding to node $v$ is $v - u + 1$ modulo $n$.
The weight of edge $(u,v)$  is given by
\begin{equation} \label{wuv.eq}
w_{uv} = \textstyle \frac{1}{2}\tr{ {M_{(1, v - u+ 1 \Mod{n})}}^{-1}},
\end{equation}

For each computation, the shift operation to obtain $M$ requires $O(n)$.
The trace of the inverse of this matrix can be found in $O(n)$ operations~\cite{RH91}. Thus, the weights of all pairs $(u,v)$ can be computed in
$O(n^3)$.  These edge weights are used in every digraph $\Gbb$ and can be looked up in constant time (for example, by storing them in an $n \times n$ matrix).

When edge weights can be computed in constant time, the digraph construction requires $O(|\Vbb| + |\Ebb|)$ operations.  Therefore, each $\Gbb$ can be constructed in $O(n^2)$.
There are $n$ such digraphs, so the construction of all digraphs requires $O(n^3)$ operations.
Finally, for each digraph, the Bellman-Ford algorithm finds the minimum-weight path of at most $k$ edges in $O(kn^2)$ time.  Thus, to find the minimum-weight path over these $n$ minimum-weight paths requires $O(kn^3)$ operations.

Combining all steps of the algorithm, we obtain a running time of $O(kn^3)$.
\end{IEEEproof}

\section{Optimal Leader Selection for Fast Convergence}~\label{algconsensus.sec}
In this section, we describe efficient algorithms that give the optimal solution to the $k$-leader selection problem for fast convergence in path and ring graphs.
Our approach is similar to that described in the previous section in that we transform the leader selection problem
into a path finding problem in a weighted digraph.  In this case, the problem is the \emph{widest path problem}~\cite{Pollack60}.
In the widest path problem, one seeks the path between two vertices for which the weight of the minimum-weight edge in that path is maximized.

\subsection{Optimal $k$-Leader Selection for a Path Graph}

Recall that for a path graph with $k$ leaders, $\Lff$ can be written in the block diagonal form in (\ref{blockL.eq}).
For a given set of leaders, the convergence rate depends on the smallest eigenvalue of $\Lff$, which in turn,
is the minimum over the smallest eigenvalues of the blocks of $\Lff$, i.e.,
\begin{align} \label{C1d.eq}
\Cerr{S} &= \min  \left[ \left\{ \lambda_{min}\left(\Lzi{u_1}\right) \right\}  \cup \left\{ \lambda_{min}\left(\Ljn{u_k}  \right)\right\} \right.  \nonumber \\
& ~~~~~~~~~~~~~~~~~~\left. \cup_{i=1}^{k-1}  \left\{  \lambda_{min}\left( \Lij{u_i}{u_{i+1}}\right) \right\}  \right].
\end{align}

As in the leader selection algorithm for coherence, we first construct a weighted digraph $\Gbc = (\Vbc, \Ebc, \Wbc)$.
The digraph has the same topology as the digraph generated in Section~\ref{algnoisepath.sec},
but the edge weights  are different.
An edge is drawn from $s$ to each node $v \in V$ with edge weight $w_{s,v} = \lambda_{min}(\Lzi{v})$.
The weight $w_{s,v}$ is the convergence rate within the subsystem consisting of nodes $1, \ldots, (v-1)$, when node $v$ is a leader
and there are no other leaders in that subsystem.
An edge is drawn from each $u \in V$ to each $v \in V$ with $v > u$ with edge weight $w_{u,v} = \lambda_{min}(\Lij{u}{v})$.
These edge weights correspond to the convergence rate within the subgraph between nodes $u$ and $v$ when both $u$ and $v$ are leaders
and no other nodes in the subgraph are leaders.
Finally, an edge is drawn from each $v \in V$ to $t$ with edge weight $w_{v,t} = \lambda_{min}(\Ljn{v})$.
The weight $w_{v,t}$ is the convergence rate within the graph consisting of nodes ${(v+1)}, \ldots, n$, when node $v$ is a leader
and there are no other leaders in that subgraph.
The weights of edges $(s,1)$, $(n,t)$, and $(u, u+1)$, $u = 1, \ldots, n-1$, are $+ \infty$.

The following proposition characterizes the relationship between the weight of a path in $\Gbc$ and $\Cerr{S}$.
The proof is similar to that of Proposition~\ref{path.prop} and is therefore omitted.
\begin{proposition} \label{widepath.prop}
Let $G$ be a path graph, and let $L$ be the respective Laplacian defined in (\ref{consensusL.eq}).
Let $\Gbc$ be the weighted digraph generated from $G$ and $L$.
Further, let $P=\{s, u_1, u_2, \ldots, u_k, t\}$ be the set of nodes on a path from $s$ to $t$ in $\Gb$ that contains $k+1$ edges
and let $w$ be the minimal edge weight in the path.  
Then, for $S  = P \setminus \{s,t\}$, $\Cerr{S} = w$.
\end{proposition}

Following from Proposition~\ref{widepath.prop},
finding the leader set for which $\Cerr{S}$ is maximized is  equivalent to finding the path from $s$ to $t$ with at most
$k+1$ edges for which the minimum edge weight is maximized.  This problem is a variation of the
widest path problem.
To find the widest path of at most $k+1$ edges efficiently, we again use a modified version of the Bellman-Ford algorithm.
Details of the modifications are given in the appendix.
The pseudocode for this algorithm is identical to that in Algorithm~\ref{optpath.alg}, except in line \ref{bf.line}, where the call
is to the modified Bellman-Ford
algorithm for the widest path rather than the minimum weight path.  We therefore omit this pseudocode for brevity.

\begin{theorem}
For a path graph $G$ with $n$ nodes, the $k$-leader selection algorithm for fast convergence identifies the leader set $S$, with $|S| \leq k$, that maximizes $\Cerr{S}$
in $O(n^3)$ operations.
\end{theorem}
\begin{IEEEproof}
The optimality of the leader set $S$ follows from Proposition~\ref{widepath.prop} and the correctness of the Bellman-Ford algorithm.
With respect to computational complexity,
the algorithm consists of two phases, each of which is performed once.  The first phase is the construction of the digraph $\Gbcr = (\Vbcr, \Ebcr)$.
with  $| \Ebcr | \in O(n^2)$. Each edge weight is given by the smallest eigenvalue of a symmetric, tridiagonal matrix.  This eigenvalue can be computed with high accuracy in $O(n)$ operations using the implicit QR algorithm of Vandebril et al.~\cite{VVM05,VVM06}.
Thus, the digraph $\Gbcr$ can be computed in $O(n^3)$ operations.

The second phase of the algorithm is the execution of the modified Bellman-Ford algorithm, which has a running time of $O(kn^2)$.
Therefore, the total running time of the algorithm is $O(n^3)$.
\end{IEEEproof}

\subsection{Optimal $k$-Leader Selection for a Ring Graph}
Using a similar approach to that for path graphs, we adapt the algorithm for optimal leader selection for coherence
to solve the $k$-leader selection problem for fast convergence in ring graphs.

To find a leader set of size at most $k$, first a single candidate node $i$ is selected as leader.  Then, the weighted digraph $\Gbbcr$
is constructed using in the same topology as described in Section~\ref{algnoisering.sec}.
The weight of each edge $(u,v)$ in $\Gbbcr$ is $w_{u,v} = \lambda_{min}(L_{(u \to v)})$.
This gives the convergence rate in the subgraph between nodes $u$ and $v$, in clockwise order, when both $u$ and $v$ are leaders
and there are no other leaders in that subgraph.

The following proposition characterizes the relationship between the weight of a path in $\Gbbcr$ and $\Cerr{S}$.
The proof is similar to that of Proposition~\ref{ring.prop} and is therefore omitted.
\begin{proposition}
Let $G$ be a ring graph, and let $L$ be the respective Laplacian as defined in (\ref{consensusL.eq}).
Let $\Gbbcr$ be the weighted digraph generated from $G$ and $L$ for a given node $i$.
Further, let $P=\{s_i, u_1, u_2, \ldots, u_{k-1}, t_i\}$ be the set of nodes on a path from $s_i$ to $t_i$ in $\Gbbcr$ that contains $k$ edges, and let $w$ be the minimal edge weight in the path.
Then, for $S  = P \setminus \{s_i,t_i\} \cup \{i\}$, $\Cerr{S} = w$.
\end{proposition}

Once the graph $\Gbbcr$ is constructed, we use our modified Bellman-Ford algorithm to find the widest path from $s_i$ to $t_i$ of at most $k$ edges.
Let $S_i$ be the set of nodes along this path.
The minimum edge weight of this path is the maximal $\Cerr{S_i}$ where $i \in S_i$ and $|S_i| \leq k$.
The optimal leader set $S^*$ is found by  finding the optimal $\Cerr{S_i}$ for all $i \in V$ and the identifying the maximum over all $i$, i.e.,
\[
S^* = \arg \max_{S_i, i \in V} \Cerr{S_i}.
\]
The pseudocode for this algorithm is nearly identical to that in Algorithm~\ref{optring.alg} and is omitted for brevity.

\begin{theorem}
For a ring graph $G$ with $n$ nodes, the  $k$-leader selection algorithm identifies the leader set $S$, with $|S| \leq k$, that maximizes $\Cerr{S}$
in $O(kn^3)$ operations.
\end{theorem}
\begin{IEEEproof}
The optimality of the leader set $S$ follows from Proposition~\ref{widepath.prop} and the correctness of the Bellman-Ford algorithm.
With respect to computational complexity, the algorithm constructs $n$ weighted digraphs, and the widest path algorithm is executed on each digraph.
As in the proof of Theorem~\ref{coherence.thm}, the $n$ digraphs can be constructed in $O(n^3)$ total operations.
The widest path Bellman-Ford algorithm is run for each digraph, requiring $O(kn^3)$ total operations.
Therefore the running time of the algorithm is $O(kn^3)$.
\end{IEEEproof}

\section{Computational Examples} \label{examples.sec}
In this section, we explore the results of our $k$-leader selection algorithms on several example graphs.
For comparison of the solution to Problem~\ref{coherence.prob}, we have implemented the greedy leader selection algorithm presented in~\cite{CABP12,CBP14}.
The greedy algorithm consists of at most $k$ iterations. In each iteration, a leader node $s$  is selected, that
when added to the leader set $S$,  yields the largest improvement in $\Merr{S}$ or $\Cerr{S}$, respectively.


For the leader selection problem for coherence, this greedy algorithm does not find the optimal leader set, but rather finds a set whose performance is within a constant factor of optimal.
Specifically, the greedy leader selection algorithm generates a  leader set $S$ of size at most $k$ such that,
\[
\Merr{S} \leq  \textstyle \left ( 1 - \left(\frac{k - 1}{k}\right)^k\right)\Mopt +  \frac{1}{e}\Mmax,
\]
where $\Mopt$ is the optimal total variance and $\Mmax \eqdef \max_{i \in V} \Merr{\{i\}}$~\cite{CBP14}.
As far as we are aware, this  algorithm is the only previously proposed solution that gives provable bounds
on the optimality of the leader set. 

We are not aware  of any previously proposed algorithms that give provable bounds on performance of the selected leader sets
for the $k$-leader selection problem for fast convergence.  
As a means for comparison, we compute the leader set using a greedy algorithm similar to that in~\cite{CABP12,CBP14}.
In each iteration $i=1, \ldots, k$, the greedy algorithm identifies the agent $s$ that yields the most increase to $\Cerr{S}$,
i.e., 
\[
s = \textstyle \arg \max_{v \in V \setminus S} \Cerr{S \cup \{v\}},
\] 
and this agent is added to the leader set.

We have implemented all algorithms in Matlab.

\subsection{Formation Coherence}
We first investigate the performance of the two leader selection algorithms for network coherence.
We consider two edge weight selection policies.  In the first, the uniform random weight policy $\nu_{ij}$ is distributed uniformly at random over 
the interval $(0.01,1)$.
In the second policy, the skewed policy, we use $\nu_{ij} = 1$ for the edges adjacent to the first half of the nodes in the ring or path
and use $\nu_{ij} = 0.01$ for the remaining edges.

\begin{figure}
\centering
  \begin{subfigure}{.9\linewidth}
          \includegraphics[scale=.38]{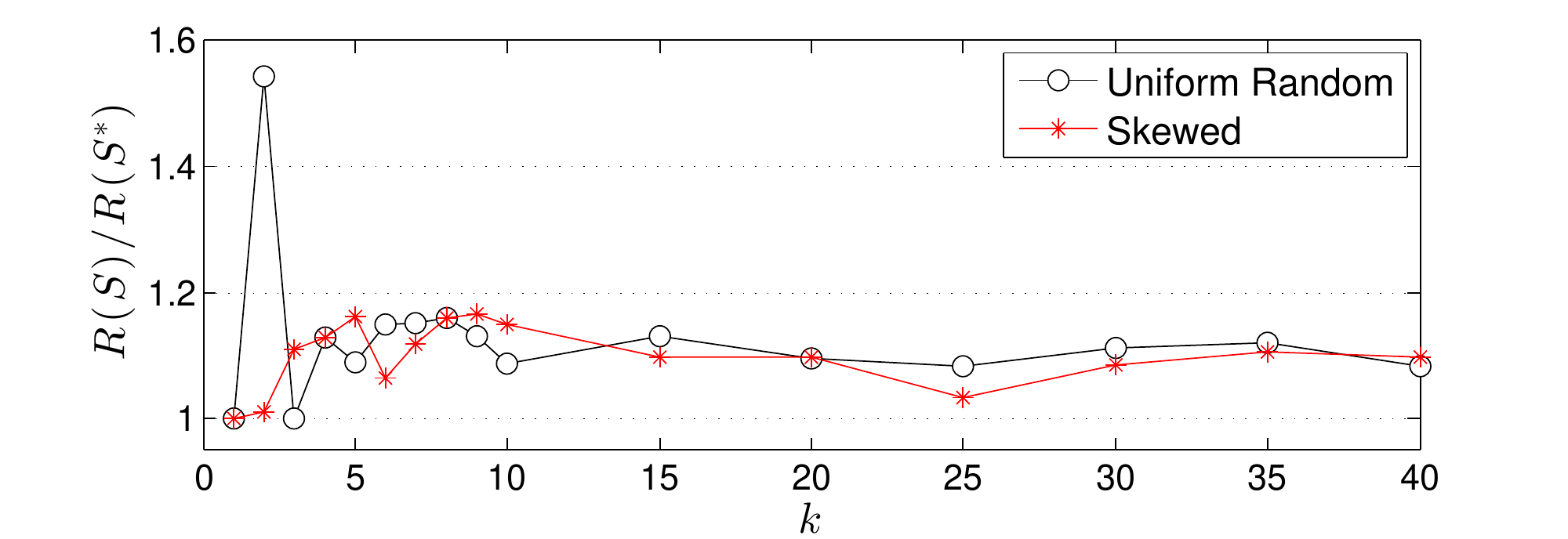}
                \caption{Total variance for leader set $S$ produced by greedy algorithm relative to the optimal leader set $S^*$
                for the uniform random and skewed edge weight policies.}\label{linerel.fig}
        \end{subfigure}
                \\
        ~ 
        \begin{subfigure}{.9\linewidth}
            \includegraphics[scale=.38]{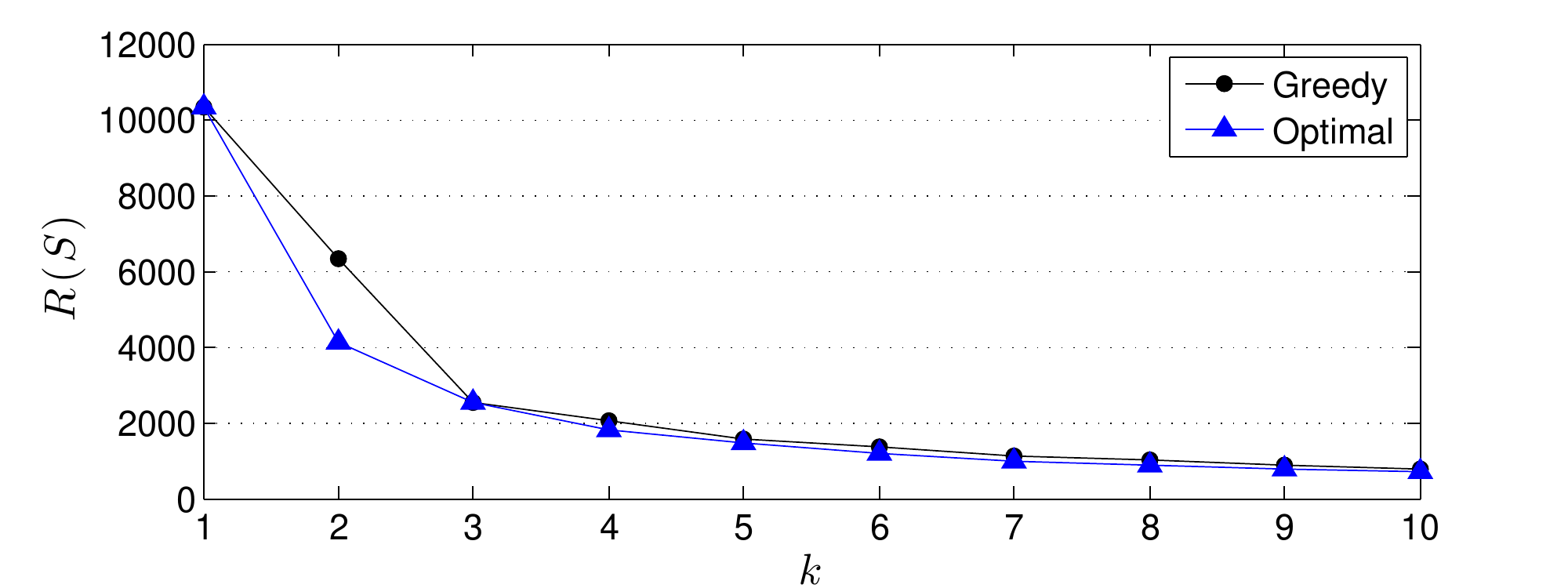}
                \caption{Total steady-state variance for the uniform random edge weight policy.}
                \label{lineabs.fig}
        \end{subfigure}%
      \caption{Comparison of formation coherence for leader sets of size $k$ on a 400 node path graph.} \label{line.fig}
\end{figure}

\begin{figure}
\centering
  \begin{subfigure}{.9\linewidth}
              \includegraphics[scale=.38]{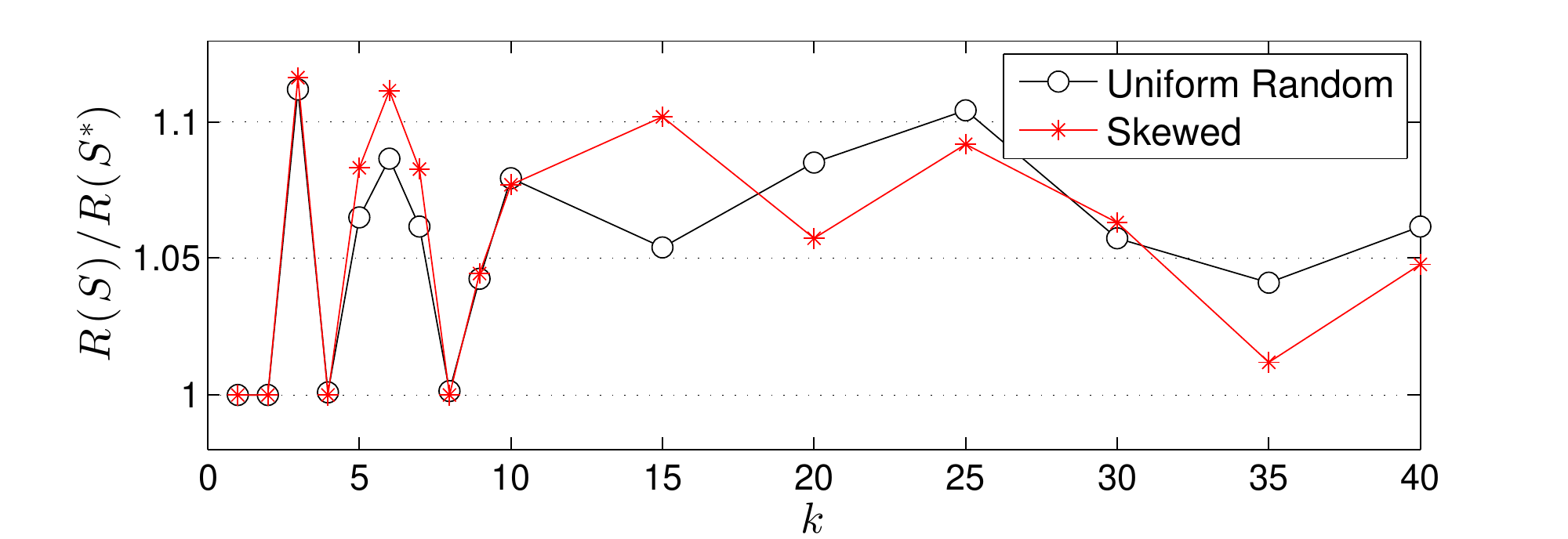}
                 \caption{Total steady-state variance for leader set $S$  produced by greedy algorithm relative to the  optimal leader set $S^*$  for the uniform random and skewed edge weight policies.}\label{ringrel.fig}
        \end{subfigure}
        \\
        ~ 
        \begin{subfigure}{.9\linewidth}
            \includegraphics[scale=.38]{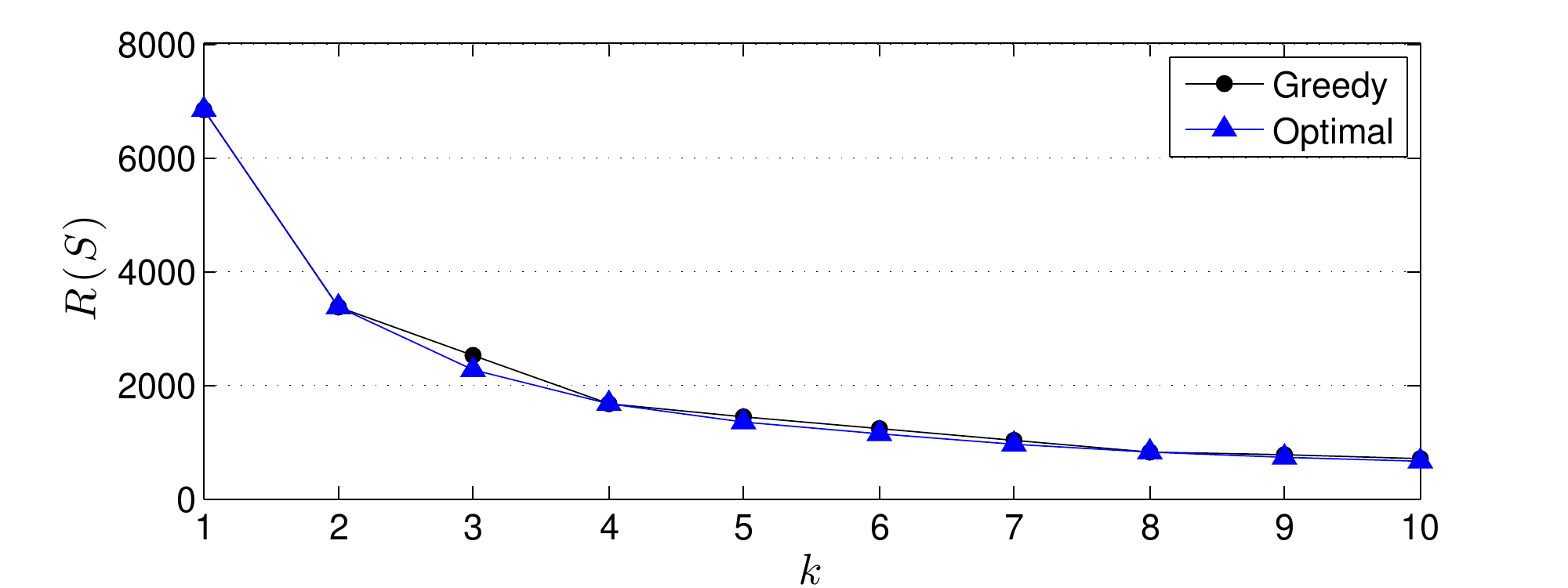}
                \caption{Total steady-state variance for the uniform random edge weight policy.}
                \label{ringabs.fig}
        \end{subfigure}%

              \caption{Comparison of formation coherence for leader sets of size $k$ on a 400 node ring graph.}\label{ring.fig}
\end{figure}

\begin{figure}
\centering
            \includegraphics[scale=.24]{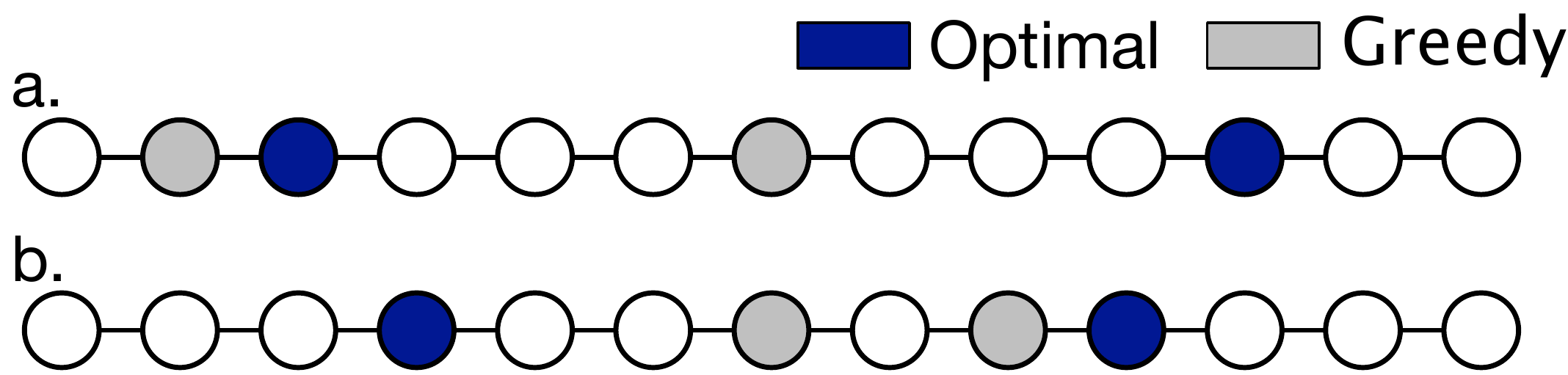}
        \caption{Leaders selected by the optimal and greedy leader selection algorithms for $k=2$ in a 13 node path graph with all edge weights equal to 1
        for (a) coherence and (b) fast convergence.}  \label{small_graph.fig}
\end{figure}

The results for a 400 node path graph are shown in Figure~\ref{line.fig}.
Figure~\ref{linerel.fig} shows the total variance $R(S)$ of the leader set selected by the greedy algorithm relative to $R(S^*)$, where $S^*$ is the optimal leader set, as found by our algorithm.
This figure shows results for both the uniform random and skewed edge weight policies.  
Figure~\ref{lineabs.fig} gives $R(S)$ for the optimal leader set and for the leader set found by the greedy algorithm where edge weights are determined using the uniform random policy.
For ${k=1}$, the optimal leader $\ell_1$ is the weighted median of the path graph (see \cite{P14}), and both algorithms select this leader.
 An interesting observation is that for ${k=2}$, the greedy algorithm demonstrates its worst relative performance for the uniform random policy.
 A reason for this can be observed in the example in Figure~\ref{small_graph.fig}a, where we show the
 optimal leaders and those selected by the greedy algorithm for ${k=2}$ in a 13 node path graph with edge weights all equal to 1.  
For the coherence problem, the locations of the optimal two leaders are symmetric.  The greedy algorithm selects the best single leader, the node in the center of the path, in the first iteration. 
It selects a node nearer to the edge of the graph in the second iteration.  The center node is a poor choice for ${k=2}$, resulting in a significantly larger total variance than the optimal.

We note that, overall, the greedy algorithm yields leader sets whose performance is fairly close to optimal.  As the number of leaders increases,
the total variance decreases for both leader selection algorithms.  The relative error of the greedy algorithm does not appear to vanish as $k$ increases.

The results for a 400 node ring graph are shown in Figure~\ref{ring.fig}.
As before, Figure~\ref{ringrel.fig} shows  $R(S)$ of the leader set selected by the greedy algorithm relative to $R(S^*)$, where $S^*$ is the optimal leader set as identified by our algorithm.
This figure shows results for both the uniform random and skewed policies.  
Figure~\ref{ringabs.fig}  shows $R(S)$ for the greedy algorithm and our optimal algorithm where edge weights are chosen using the uniform random policy.
For  all $k$ greater than 1, the greedy algorithm selects a sub-optimal leader set. 
For larger values of $k$,
the performance of the greedy algorithm appears to stabilize around 1.05 times the optimal $R(S^*)$ for both edge weight policies.
We note that for larger $k$, the performance of the greedy algorithm is similar in the ring and path graphs.

\subsection{Fast Convergence}
We next explore the impact of leader set selection  on convergence rate.
We use a  policy where $W_{ij}$ is drawn uniformly at random from $(0,1)$ and a skewed policy where $W_{ij} = 1$ for edges adjacent to 
the first half of the follower nodes in the path or ring and $W_{ij} = 100$ for the remaining edges.

\begin{figure}
\centering
    \begin{subfigure}{.9\linewidth}\centering
                   \includegraphics[scale=.38]{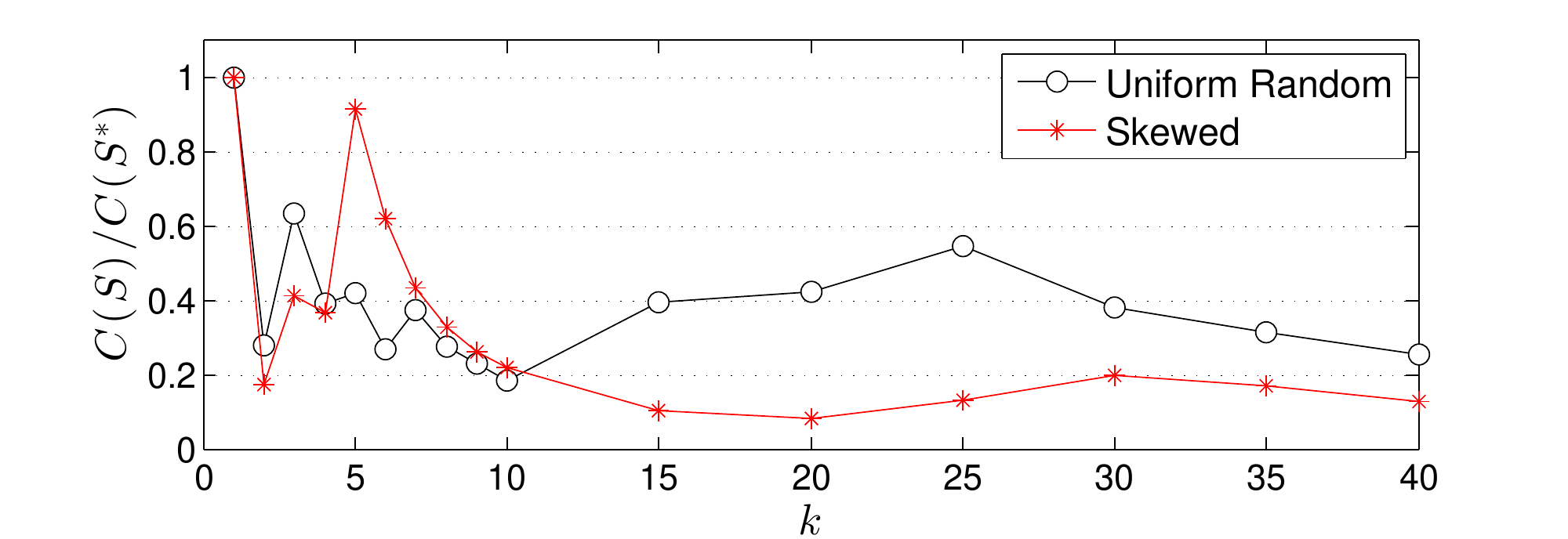}
                \caption{Convergence rate for leader set $S$ produced by the greedy algorithm relative to  the optimal leader set $S^*$
            for the uniform random and skewed edge weight policies.}\label{lineraterel.fig}
        \end{subfigure} 
              \\      
                ~ 
    
        \begin{subfigure}{.9\linewidth}\centering
            \includegraphics[scale=.38]{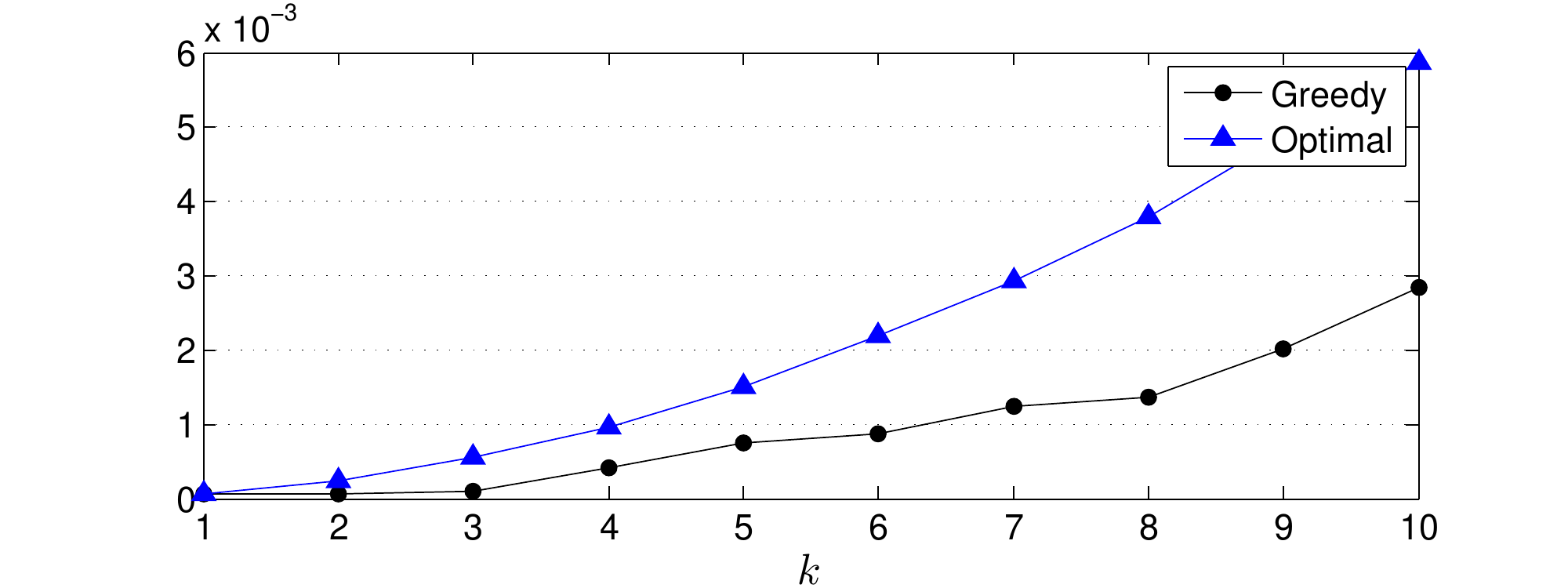}
                \caption{Convergence rate for the uniform random edge weight policy.}
                \label{linerateabs.fig}
        \end{subfigure}%
  
        \caption{Comparison of convergence rate for leader sets of size $k$  on a 400 node path graph.}\label{linerate.fig}

    \end{figure}

The results for a 400 node path graph are shown Figure~\ref{linerate.fig}.
Figures~\ref{lineraterel.fig} shows $\Cerr{S}$, i.e., the value of the minimal eigenvalue of $\Lff$,  for the leader set selected by the greedy algorithm relative to $\Cerr{S^*}$, where $S^*$ is the optimal leader set resulting from our algorithm.   Results are shown for both edge weight policies.
Figure~\ref{linerateabs.fig} shows $\Cerr{S}$ for the optimal leader set
and leader set identified by the greedy algorithm for various leader set sizes $k$ under the uniform random edge weight policy. 
A larger $C(S)$ corresponds to faster convergence.
As in the previous set of experiments, both algorithms find the optimal leader for $k=1$.   
For $k$ greater than 1, the greedy algorithm performs poorly in most cases.  Further, the relative error of the greedy algorithm does not appear to vanish as $k$ increases.

\begin{figure}
\centering
    \begin{subfigure}{.9\linewidth}\centering
                   \includegraphics[scale=.37]{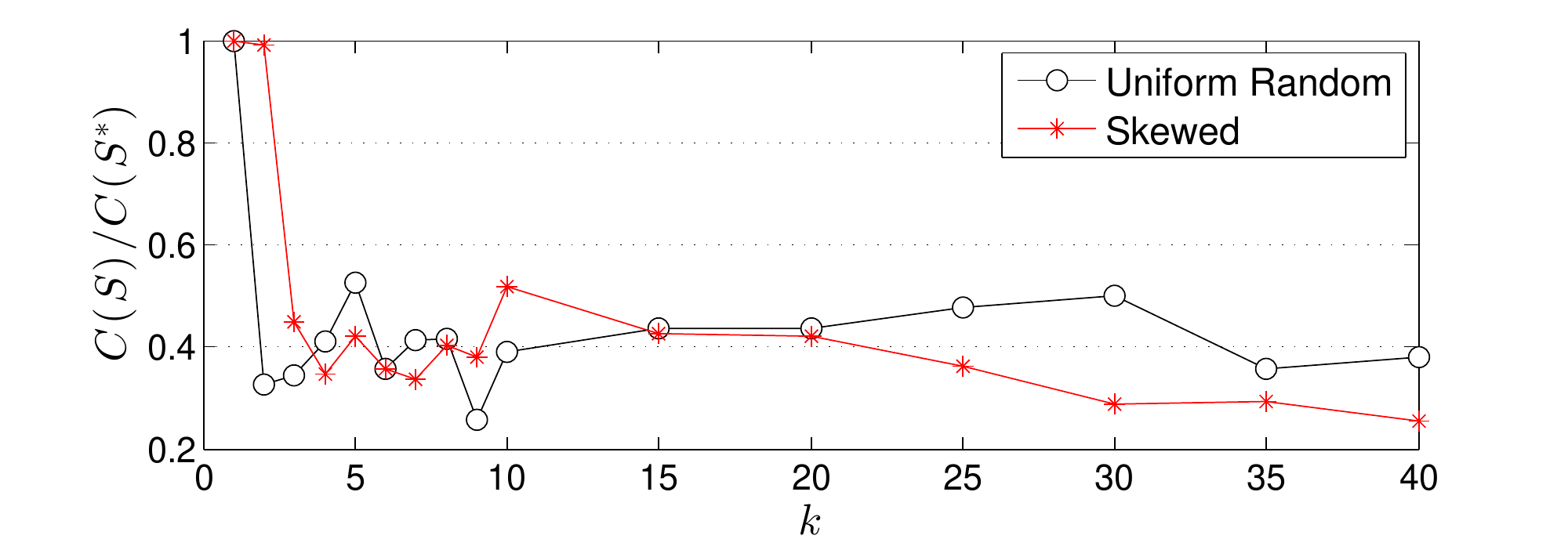}
                \caption{Convergence rate for leader set $S$ produced by greedy algorithm relative to the optimal leader set $S^*$
                for the uniform random  and skewed edge weight policies.}\label{ringraterel.fig}
        \end{subfigure} 
              \\      
                ~ 
    
        \begin{subfigure}{.9\linewidth}\centering
            \includegraphics[scale=.38]{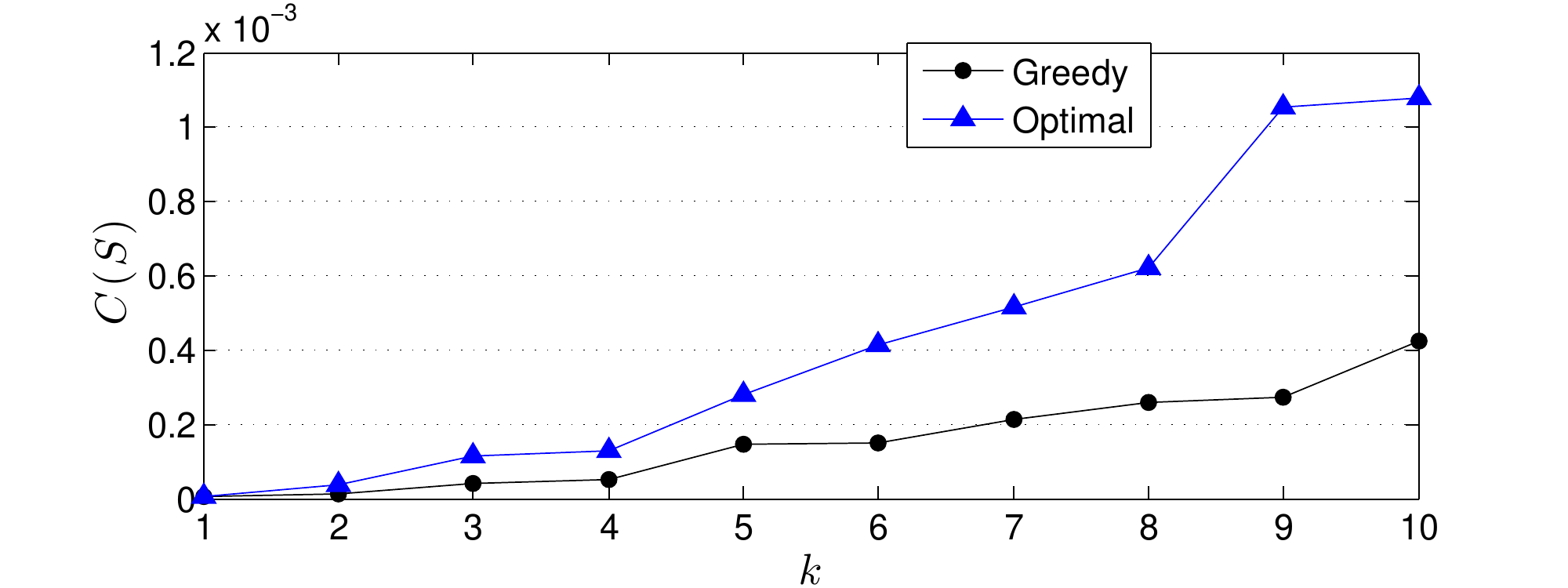}
                \caption{Convergence rate for the uniform random edge weight policy.}
                \label{ringrateabs.fig}
        \end{subfigure}%
  
      \caption{Comparison of convergence rate for leader sets of size $k$  on a 400 ring path graph.}\label{ringrate.fig}
\end{figure}

The results for a 400 node ring graph are shown in Figure~\ref{ringrate.fig}.
As above, Figure~\ref{ringraterel.fig} shows  $\Cerr{S}$ of the leader set selected by the greedy algorithm relative to the optimal $\Cerr{S^*}$
for both edge weight policies.
Figure~\ref{ringrateabs.fig} gives $\Cerr{S}$ for the greedy algorithm and optimal leader set for various leader set sizes $k$ under the uniform random edge weight policy.
For both path and ring graphs, the greedy algorithm is less effective in solving the $k$-leader selection for fast convergence than for solving
the $k$-leader selection problem for coherence.

Finally we note that the optimal leader set for coherence is not necessarily the same as the optimal leader set for fast convergence.  This can be observed
in the 13 node network shown in Figure~\ref{small_graph.fig}, where edge weights are all equal to 1.

\section{Conclusion} \label{conclusion.sec}
We have investigated the problem of optimal $k$-leader selection in leader-follower consensus systems
for two performance objectives, coherence and convergence rate. 
A na\"{\i}ve solution to either leader selection problem has combinatorial complexity, however, it is unknown whether these problems are
NP-Hard or if efficient polynomial-time solutions can be found.
In this work, we have taken a step towards addressing this open question.
We have shown that, in one-dimensional undirected, weighted graphs, namely path graphs and ring graphs, both the $k$-leader selection problem for coherence
and the $k$-leader selection problem for fast convergence can be solved in polynomial time in both $k$ and the network size $n$.
Further, for each problem, we have given an $O(n^3)$ solution for optimal $k$-leader selection in path graphs and  an $O(kn^3)$ solution for optimal $k$-leader selection in ring graphs.

While our approach depends on the specific structure of path and ring graphs and thus cannot be easily extended to other graph topologies,
we anticipate that by applying other techniques used for network facility location, it will be possible to develop efficient algorithms for additional topologies such as tree graphs.
We plan to address this in future work. Finally, we will explore
using similar algorithmic techniques for leader selection in other dynamics, including networks where leaders are also subject to stochastic noise.

\section*{Appendix} 

Pseudocode for our modified Bellman-Ford algorithm for finding a minimum weight path of at most $H$ edges in a digraph is given in Algorithm~\ref{mbf.alg}.
The algorithm performs $H$ iterations; in each iteration $m$, it finds the minimum-weight path of exactly $m$ hops from node $s$ to every other node.
If no such path exists, the path weight is infinite. 
\begin{algorithm}[t]

\caption{Modified Bellman-Ford algorithm.}\label{mbf.alg}
\begin{algorithmic}[1]\footnotesize
\State \textbf{Input:} weighted digraph ${\cal G} = ({\cal V},{\cal E},{\cal W})$, source node $s$,  target node $t$,
number of edges $H$
\State \textbf{Output:} List $P$ of vertices on min. weight path from $s$ to $t$ of at most $H$ edges, weight $w$ of path $P$
\Statex~
\Statex \emph{ /* $d_{m}(v)$ stores the weight of the minimum-weight path from node $s$ to node $v$ that contains exactly $m$ edges. */}
\State $d_{0}(s) \gets 0$ \label{d0init.line}
\State $d_{m}(s) \gets \infty$ for $m = 1\ldots H$ \label{dsinit.line}
\State $d_{m}(v) \gets \infty$ for $v \in V, v \neq s, m = 0 \ldots H$ \label{dmvinit.line}
\Statex~
\State \emph{ /* $\pi_m(v)$ stores the ID of the predecessor of node $v$ in the minimum-weight path of $m$ edges. */ }
\State $\pi_{m}(v) \gets \perp$ for $v \in V, m = 0 . . . H$
\Statex~
\Statex \emph{/* Path $P$ initialized to empty list */}
\State $P \gets [~]$
\Statex~
\Statex \emph{/*~Find minimum-weight paths of all lengths up to $H$ from $s$ to all other nodes.~*/}
\For{$m = 1 \ldots H$}
\For{$v \in V$}
\State $d_{m}(v) \gets \displaystyle \min_{u \in N(v)}\{d_{m-1}(u) + w(u,v) \}$ \label{dupdate.line}
\State $\pi_{m}(v) \gets  \displaystyle  \arg\min_{u \in N(v)}\{d_{m-1}(u) + w(u,v) \}$  \label{piupdate.line}
\EndFor
\EndFor
\Statex~
\Statex \emph{/* Identify min. weight path of at most $H$ edges. */}
\State $ L = H$
\While{$d_{L}(t) = d_{L-1}(t)$}
	\State $L \gets L-1$
\EndWhile
\State $P$.append($t$);
\State $u \gets \pi_L(t)$
\State $P$.append($u$)
\For{$i=(L-1) \ldots 2$}
	\State $u \gets \pi_{i}(u)$
	\State $P$.append($u$)
\EndFor
\State \textbf{return} $(P,d_L(t))$
\end{algorithmic}
\end{algorithm}

To find the path of length $m$ from  $s$ to a node $v$, the algorithm examines each node $u \in N(v)$ and finds minimum over the weight of path from $s$ to $u$ of length ${m-1}$ plus the weight of edge $(u,v)$.   Each iteration of the algorithm requires a number of operations of the order of the number of edges in the graph.  After $H$ iterations, the algorithm finds the minimum-weight path from $s$ to $t$ of at most $H$ edges over the computed paths from $s$ to $t$ of $m$ edges, $m = 1\dots H$.

With small modifications, Algorithm~\ref{mbf.alg} can be used the find the widest path in the digraph ${\cal G}$ of at most $H$ edges.
To solve the widest path problem, the algorithm keeps track of the minimum edge weight along each path instead of the total weight of the path.
Specifically,  the initialization of the path weights in lines \ref{d0init.line} - \ref{dmvinit.line} is changed to:
\begin{align*}
& d_{0}(s) \gets \infty \\
& d_{m}(s) \gets -\infty~~\text{for}~m = 1\ldots H \\
& d_{m}(v) \gets -\infty~~\text{for}~v \in V, v \neq s, m = 0 \ldots H,
\end{align*}
and the update of the path weights in lines \ref{dupdate.line} and \ref{piupdate.line} becomes:
\begin{align*}
 & d_{m}(v)\gets \textstyle \max_{u \in N(v)}\{ \min  \{d_{m-1}(u),w(u,v)\} \} \\
& \pi_{m}(v) \gets  \textstyle  \arg\max_{u \in N(v)}\{ \min \{d_{m-1}(u), w(u,v)\} \}.
\end{align*}

The running time of Algorithm~\ref{mbf.alg}, for both the minimum-weight path problem and the widest path problem, is $O(H |{\cal E}|)$.  We refer the reader to~\cite{CLRS10} for more details on the original Bellman-Ford algorithm and its complexity analysis.

%
\bibliographystyle{IEEETran}
\bibliography{leader}

\end{document}